\documentclass[conference]{IEEEtran}
\IEEEoverridecommandlockouts

\usepackage{graphicx}
\graphicspath{{./images/}}
\usepackage[caption=false]{subfig}
\usepackage{amsmath, xparse}
\usepackage{mathtools}
\mathtoolsset{showonlyrefs}
\usepackage[dvipsnames]{xcolor}
\usepackage{amssymb}
\usepackage{amsthm}
\usepackage{breqn}
\usepackage{enumitem}

\newtheorem{theorem}{Theorem}[section]
\newtheorem{assumption}{Assumption}[section]
\newtheorem{definition}{Definition}[section]
\newtheorem{proposition}{Proposition}[section]
\newtheorem{remark}{Remark}[section]
\newtheorem{example}{Example}[section]

\title{\LARGE \bf $\alpha$-stability of Differentially Flat Systems with Application to Newton-Raphson Tracking Control for Vehicle Dynamics}

\author{Aadila Ali Sabry and Gennaro Notomista%
	\thanks{This work has been partially supported by the NSERC Alliance Grant ALLRP 599163 - 24.}%
	\thanks{The authors are with the Department of Electrical and Computer Engineering, University of Waterloo, Waterloo, ON, Canada {\tt\footnotesize \{aalisabr@uwaterloo.ca,gennaro.notomista\}@uwaterloo.ca}.}%
}

\begin{document}

\maketitle
\thispagestyle{empty}
\pagestyle{empty}

\begin{abstract}

This paper studies the $\alpha$-stability property of differentially flat nonlinear dynamical systems. The results build off the recently introduced notion of $\alpha$-stability, which is particularly amenable to characterize the ability of a system to track dynamic output reference signals. We consider systems controlled using the Newton-Raphson tracking controller, which results in closed-form control policies, therefore it is computationally efficient, and it has been shown to be effective to control a large variety of mobile robots, including autonomous vehicles. The main results of the paper consist in sufficient conditions for the $\alpha$-stability of differentially flat systems and for the equivalence between the proposed control algorithm and the Newton-Raphson tracking controller applied directly to the nonlinear dynamics. We demonstrate the behavior of the proposed controller applied to the kinematic unicycle and dynamic bicycle models.

\end{abstract}

\section{Introduction}

Desirable properties of tracking controllers for autonomous systems are robust performance guarantees and low computational complexity. This is particularly true in safety-critical applications, such as collision mitigation in autonomous driving, where the decision making process has to take place in short amounts of time, yet safety guarantees have to hold for all times \cite{wang2020safety,hicks2018safety}. While many tracking controllers have been developed and successfully applied in real-world scenarios, the control systems and transportation community are still interested in improving closed-loop controlled systems performance while minimizing the computational resources required to compute the controller.

Model predictive control is a prominent example of a general purpose framework for constraint satisfaction where significant effort has been put in reducing computational complexity, based, e.g., on successive linearizations of the system dynamics \cite{allgower2012nonlinear,borrelli2017predictive}. The Newton-Raphson (NR) controller \cite{wardi2024tracking} has been recently developed in order to guarantee robust tracking for generic nonlinear systems, leveraging closed-form integral control actions. The NR controller has been shown to be successful on a variety of cyber-physical systems applications, including autonomous vehicles \cite{notomista2024safe}, flying robots \cite{morales2024newton}, and multi-robot systems \cite{shivam2019tracking}. However, the performance guarantees for systems modeled by generic nonlinear dynamics cannot be easily verified. In particular, sufficient conditions for $\alpha$-stability---which quantifies the robust tracking performance of the NR controller---are only provided for linear system dynamics.

In this paper, we propose a strategy to apply the NR controller to a class of nonlinear systems known as differentially flat systems \cite{martin2006flat}. Differentially flat systems typically lend themselves to elegant control solutions since the control inputs can be computed algebraically from trajectories of the so-called \textit{flat outputs}. The authors in \cite{mellinger2011minimum}, for instance, are able to achieve high-performance control of dynamics as complicated as those of a quadrotor, thanks to differential flatness of the mathematical model. More general optimal control solutions have been recently proposed in \cite{beaver2024optimal}, where the authors propose a strategy to solve optimal controllers for differentially flat systems. Nevertheless, per-step optimization within the feedback loop is computationally expensive.

Our proposed approach consists of applying a NR controller to the trivial (linear) dynamics of the flat outputs of the system and leveraging the endogenous (algebraic, closed-form) transformation in order to compute the control input signals to supply to the system. Closed-loop convergence is shown for a modification of this controller, and sufficient conditions are derived for the $\alpha$-stability of the NR controller in the flat output dynamics. In some cases, this controller is equivalent to applying the NR controller directly to the nonlinear dynamics. We characterize a set of conditions under which this occurs.

The control framework we propose in this paper is particularly amenable for complicated nonlinear dynamics, such as those modeling nonholonomic mobile robots. We therefore showcase the behavior of our controller first, for illustration purposes, on a kinematic unicycle model, and then on a dynamic bicycle model. The latter is applicable to control applications involving autonomous vehicles and, thanks to its low computational complexity, can be employed in safety-critical scenarios where controllers need to respond in a timely fashion to mitigate collision or realize evasive maneuvers \cite{skibik2023mpc}.

The rest of the paper is organized as follows. The next section introduces the mathematical concepts required for the development of the proposed control framework and for the proof of the sufficient conditions for tracking performance of our controller. In Section~\ref{sec:main} we illustrate the main results of the paper:
\begin{enumerate}[label=(\roman*)]
	\item Sufficient conditions for $\alpha$-stability of the flat output dynamics
	\item Proof of convergence of a modified version of the designed tracking controller
	\item Characterization of the equivalence between the proposed control strategy and a NR controller applied directly to the nonlinear dynamics under certain assumptions
\end{enumerate}
Section~\ref{sec:examples} illustrates how to apply the proposed controller to kinematic unicycle and dynamic bicycle motion models, while Section~\ref{sec:simulations} reports the results of simulations performed using these two models. Section~\ref{sec:conclusions} concludes the paper.

\section{Mathematical Background}
\label{sec:background}

\subsection{Newton-Raphson Controller} \label{sec:background:NR}

\begin{figure}
\centering
\includegraphics[width=0.5\textwidth]{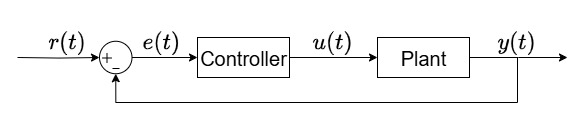}
\caption{Controlled System}
\label{fig:system}
\end{figure}

Consider a system controlled using the typical feedback control paradigm, as depicted in Fig.~\ref{fig:system}, where the input signal $u(t)$, reference $r(t)$ and output $y(t)$ are vectors in $\mathbb{R}^m$ for all times $t$, and the state $x(t) \in \mathbb{R}^n$. The plant evolves according to the generic nonlinear dynamics:
\begin{align}
	\begin{cases}
		\dot{x}(t) = f(x(t), u(t))\\
		y(t) = h(x(t)),
	\end{cases}
    \label{eqn:plant}
\end{align}
where $f$ and $h$ are functions of domains and codomains of appropriate dimensions. The NR controller we employ in this paper has been introduced in \cite{wardi2024tracking} and it has been derived based on the Newton-Raphson method of iterative root finding to the equation $r (t + T) - \hat{y} (t + T) = 0$, where $\hat{y}(t+T) = g(x(t), u(t))$ is a \textit{prediction} of the output at time $t+T$, made at time $t$. $T$ is typically referred to as the prediction horizon. The dynamically-defined NR controller is defined as:
\begin{align}
    \dot{u}(t) &= \alpha \left( \frac{\partial g}{\partial u} \right)^{-1} e(t+T), \label{eqn:controller}
\end{align}
where $e(t+T) := r(t+T) - \hat{y}(t+T)$, and $\alpha > 1$. The latter is referred to as the \textit{speedup factor} and it has been shown that, for certain example systems, increasing $\alpha$ while keeping $T$ fixed can stabilize the system. Therefore, a suitable notion of stability should take into account the behavior of the system as $\alpha \rightarrow \infty$.

To define this notion of stability, some notation needs to be introduced. Consider the closed-loop system described by the plant \eqref{eqn:plant} controlled by \eqref{eqn:controller}, which can be viewed as having an augmented state $z(t) := \begin{bmatrix} x(t)^\top & u(t)^\top \end{bmatrix} ^\top \in \mathbb{R}^{n+m}$. Denote by $L^\infty(\mathbb{R}^m)$ the space of essentially-bounded functions $p : [0, \infty) \rightarrow \mathbb{R}^m$ and by $C$ the space of bounded, continuous functions $r : [0, \infty) \rightarrow \mathbb{R}^m$ having a bounded, piecewise-continuous derivative. For a given initial condition $z_0 \in \mathbb{R}^{n+m}$, the state trajectory $\{z(t): t \in [0, \infty)\}$ is \textit{well-defined} if it comprises the unique continuous, piecewise-continuously differentiable solution of \eqref{eqn:plant}--\eqref{eqn:controller} over $t \in [0, \infty)$ with $z(0) = z_0$. A well-defined state trajectory is \textit{nonsingular} if for every $z(t) = \begin{bmatrix} x(t)^\top & u(t)^\top \end{bmatrix} ^\top$, the partial Jacobian $\frac{\partial g}{\partial u}(x(t), u(t))$ is nonsingular. The following definition, taken from \cite{wardi2024tracking}, is required in the remainder of the paper in order to talk about $\alpha$-stability and prove sufficient conditions for it to hold in the case of differentially flat systems.

\begin{definition}[Definition 1 in \cite{wardi2024tracking}]
    Given an open set $\Omega \subset L^\infty(\mathbb{R}^m)$ and a closed set $\Gamma \subset \mathbb{R}^{n+m}$, the system defined by \eqref{eqn:plant} and \eqref{eqn:controller} is $\alpha$-\textit{stable} with respect to $\Omega$ and $\Gamma$ if there exist $\overline{\alpha} > 0$ and two continuous, monotone-nondecreasing functions $\gamma: \mathbb{R}^+ \rightarrow \mathbb{R}^+$ and $\beta: \mathbb{R}^+ \rightarrow \mathbb{R}^+$ such that for every $\alpha > \overline{\alpha}$, reference $\{r(t)\} \in \Omega \cap C$, and initial condition $z_0 := z(0) \in \Gamma$, the trajectory $\{z(t)\}$ is well-defined and nonsingular, and the following condition is satisfied:
    \begin{align}
        \| z \|_\infty \leq \gamma(\| r \|_\infty) + \beta(\| z_0 \|).
    \end{align}
\end{definition}
Under mild assumptions on $f$, $g$, $h$ and $r$, $\alpha$-stability implies tracking convergence of the controlled system \eqref{eqn:plant}--\eqref{eqn:controller}.

The following assumptions guarantee the existence of a unique continuous piecewise-differentiable solution for \eqref{eqn:plant} on the time horizon $\{t : t \geq 0\}$ as long as the input $u(t)$ is piecewise continuous and bounded, and the closed-loop system \eqref{eqn:plant}--\eqref{eqn:controller} is well-posed.
\begin{assumption}[Assumption 1 in \cite{wardi2024tracking}] \label{assumption1}
\par\noindent
    \begin{enumerate}
        \item The function $f$ is continuously differentiable, and for every compact set $\Gamma \subset \mathbb{R}^m$, there exists $K > 0$ such that for $(x, u) \in \mathbb{R}^n \times \Gamma$,
        \begin{align}
            \| f(x, u) \| \leq K (\| x \| + 1).
        \end{align}
        \item The function $h$ is continuously differentiable.
    \end{enumerate}
\end{assumption}
\begin{assumption}[Assumption 2 in \cite{wardi2024tracking}] \label{assumption2}
    The function $g$ is continuously differentiable in $(x, u)$, and its partial Jacobian, $\frac{\partial g}{\partial u}(x, u)$ is locally Lipschitz continuous in $(x, u)$.
\end{assumption}

The following are the main results from \cite{wardi2024tracking} on the tracking performance guarantees of a system controlled using the NR controller.
\begin{proposition}[Proposition 3 in \cite{wardi2024tracking}]
    Consider the closed-loop system defined by \eqref{eqn:plant}--\eqref{eqn:controller}. Suppose that assumptions \eqref{assumption1} and \eqref{assumption2} are satisfied. Given an open set $\Omega \subset L^\infty(\mathbb{R}^m)$ and a closed set $\Gamma \subset \mathbb{R}^{n+m}$, suppose that the system is $\alpha$-stable with respect to $\Omega$ and $\Gamma$. Then for every exogenous reference input $\{r(t)\} \in \Omega \cap C$, and for every $z_0 \in \Gamma$,
    \begin{align}
        \lim_{\alpha \to \infty} \limsup_{t \to \infty} \| r(t) - \hat{y}(t) \| = 0
    \end{align}
\end{proposition}
In this paper we will take the predictor to be of the following form: let $\{ \zeta(\tau) : \tau \in [t, t+T] \}$ be defined by the differential equation $\dot{\zeta}(\tau) = f(\zeta(\tau), u(t))$ and the boundary condition $\zeta(t) = x(t)$. Define $\hat{y}(t+T) = h(\zeta(t+T))$. For a linear system with dynamics
\begin{align}
    \dot{x}(t) &= Ax(t) + Bu(t), \: y(t) = Cx(t) \label{eq:linear-system-dynamics}
\end{align}
the output predictor function takes the form $\hat{y}(t+T) = Ce^{AT}x(t) + CA^{-1} \int_0^T e^{A(T-\tau)} d \tau B u(t)$. The control input is then obtained by evaluating the following expression:
\begin{align}
    \dot{u}(t) = \alpha \left( C \int_0^T e^{A \tau} d \tau B \right)^{-1} \Bigg( &r(t+T) - Ce^{AT} x(t) \\
    &- C \int_0^T e^{A \tau} d \tau Bu(t) \Bigg),
\end{align}
and the augmented state $\begin{bmatrix}x(t)^\top & u(t)^\top\end{bmatrix}^\top$ evolves according to the following linear dynamics:
\begin{align}
    \begin{bmatrix}
        \dot{x}(t) \\
        \dot{u}(t) \\
    \end{bmatrix} &= \begin{bmatrix}
        A & B \\
        - \alpha \left( C \int_0^T e^{A \tau} d \tau B \right)^{-1} C e^{AT} & - \alpha I \\
    \end{bmatrix} \begin{bmatrix}
        x(t) \\
        u(t) \\
    \end{bmatrix} \\
    &+ \begin{bmatrix}
        0 \\
        \alpha \left( C \int_0^T e^{A \tau} d \tau B \right)^{-1} \\
    \end{bmatrix} r(t+T).
    \label{eqn:augmented_state_evolution}
\end{align}

Denote by $P_\alpha(s)$ the characteristic polynomial of the matrix coefficient of $\begin{bmatrix}x (t)^\top & u(t)^\top \end{bmatrix} ^\top$ in Eq. \eqref{eqn:augmented_state_evolution}, and let $P_i(S)$ be the coefficient of $\alpha^{m-i}$ in $P_\alpha(s)$; $P_i(s)$ is a polynomial in $s$. Let $\tilde{P}_i(s)$ be the highest degree term in $P_i(s)$, and define $Q(S) := \frac{1}{s^n}\sum_{i=0}^m \tilde{P}_i(s)$. Note that $P_0(s), Q(s)$ are independent of $\alpha$. We conclude this section by recalling the following sufficient condition for $\alpha$-stability that holds for linear systems.
\begin{theorem}[Theorem 1 in \cite{wardi2024tracking}] \label{theorem:suff-conds}
    Consider the system defined by Eq. \eqref{eqn:augmented_state_evolution}. If the polynomials $P_0(s)$ and $Q(s)$ have all of their roots in the open left-half plane (LHP), then the system is $\alpha$-stable with respect to $\Omega := L^\infty(\mathbb{R}^m)$ and $\Gamma := \mathbb{R}^{n+m}$.
\end{theorem}

\subsection{Differential Flatness} \label{sec:background:flatness}

\begin{figure*}
	\centering
	\subfloat[]{%
		\label{fig:a}
		\includegraphics[width=0.48\linewidth]{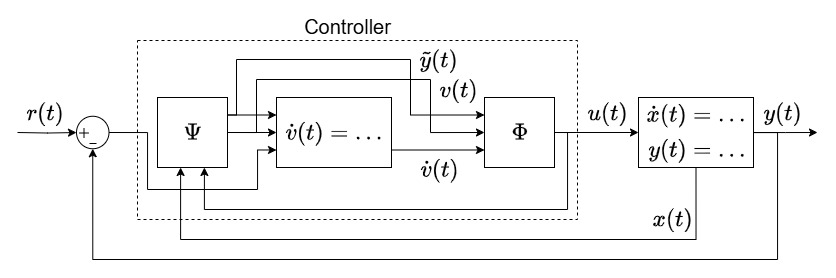}%
	}\hfill
	\subfloat[]{%
		\label{fig:b}
		\includegraphics[width=0.48\linewidth]{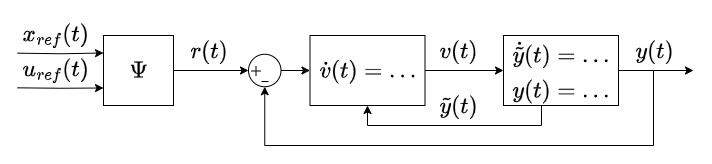}%
	}
	\caption{(a) Block diagram depicting the sequence (1) endogenous transformation to flat dynamics, (2) NR in flat dynamics and (3) inverse endogenous transformation to real dynamics, as described in Section~\ref{subsec:closed-loop-convergence}. $\Psi$ denotes the endogenous transform and $\Phi$ its inverse. (b) Block diagram with flat dynamics evolving in parallel with the real dynamics; the trajectory generated for the real dynamics are passed through the endogenous transformation to generate the reference in the flat dynamics (Section~\ref{subsec:NR-on-original-dynamics}), with $\Psi, \Phi$ as in Figure \ref{fig:a}.}
	\label{fig:placeholder}
\end{figure*}

To introduce differential flatness, we will first define systems over infinite-dimensional spaces. We will work over the \textit{Fr\'echet topology}: a vector field in this topology is \textit{smooth} if each component of the vector field depends on a finite but arbitrary number of variables and is smooth in the usual sense.

\begin{definition}[Definition 1 in \cite{martin2006flat}]
    A \textit{system} is a pair $(\mathfrak{M}, F)$ where $\mathfrak{M}$ is a smooth manifold (possibly of infinite dimension), and $F$ is a smooth vector field on $\mathfrak{M}$. A \textit{trajectory} of the system is a mapping $t \mapsto \xi(t)$, where $\xi \in \mathfrak{M}$ such that $\dot{\xi}(t) = F(\xi(t))$.
\end{definition}

Under this definition, the system in \eqref{eqn:plant} is defined by the pair $(\mathbb{R}^n \times \mathbb{R}^m \times \mathbb{R}^\infty_m, F)$, where $\mathbb{R}_m^\infty = \mathbb{R}^m \times \mathbb{R}^m \times \dots$ denotes the product of a countably infinite number of copies of $\mathbb{R}^m$ and $F(x(t), u(t), \dot{u}(t), \dots) = (f(x(t), u(t)), \dot{u}(t), \ddot{u}(t), \dots)$. This definition makes the concept of state dimension meaningless: in particular, a system and its dynamic extension have the same description. The following example will be useful for the development of the proposed controller, which hinges on the assignment of trivial (chain of integrators) dynamics to the dynamics of the flat outputs of a differentially flat system.

\begin{example}[Example 1 from \cite{martin2006flat}]
    The \textit{trivial system} $(\mathbb{R}_m^\infty, F_m)$ with coordinates $(y, \dot{y}, \ddot{y}, \dots)$ and vector field $F_m(y, \dot{y}, \ddot{y}, \dots) = (\dot{y}, \ddot{y}, \dddot{y}, \dots)$ describes any system comprised of $m$ chains of integrators of arbitrary lengths.
\end{example}

Consider two systems $(\mathfrak{M}, F)$ and $(\mathfrak{N}, G)$ and a smooth mapping $\Psi: \mathfrak{M} \rightarrow \mathfrak{N}$. If $t \mapsto \xi(t)$ is a trajectory of $(\mathfrak{M}, F)$, the composed mapping $t \mapsto \zeta(t) = \Psi(\xi(t))$ satisfies the chain rule $\dot{\zeta}(t) = \frac{\partial \Psi}{\partial \xi}(\xi(t)) \cdot F({\xi}(t))$.

\begin{definition}[\cite{martin2006flat}] \label{def:psi-related}
    Given two systems $(\mathfrak{M}, F)$ and $(\mathfrak{N}, G)$ and a smooth mapping $\Psi: \mathfrak{M} \rightarrow \mathfrak{N}$, $F$ and $G$ are said to be $\Psi$-\textit{related} if for all $\xi$, $G(\Psi(\xi)) = \frac{\partial \Psi}{\partial \xi}(\xi) \cdot F({\xi})$. $\Psi$ is an \textit{endogenous transformation} if it has a smooth inverse.
\end{definition}

Definition \ref{def:psi-related} says that $F, G$ are $\Psi$-related if trajectories $t \mapsto \xi(t)$ of $(\mathfrak{M}, F)$ are mapped to trajectories $t \mapsto \Psi(\xi(t))$ of $(\mathfrak{N}, G)$ by $\Psi$. If $\Psi$ is endogenous with inverse $\Phi$, $F$ and $G$ are also $\Phi$-related, and there is a one-to-one correspondence between the trajectories of the two systems.

\begin{definition}[Definition 2 from \cite{martin2006flat}]
    Two systems $(\mathfrak{M}, F)$ and $(\mathfrak{N}, G)$ are \textit{equivalent at} $(p, q) \in \mathfrak{M} \times \mathfrak{N}$ if there exists an endogenous transformation from a neighborhood of $p$ to a neighborhood of $q$. $(\mathfrak{M}, F)$ and $(\mathfrak{N}, G)$ are said to be \textit{equivalent} if they are equivalent at every pair of points $(p, q)$ in a dense open subset of $\mathfrak{M} \times \mathfrak{N}$.
\end{definition}

We will use the shorthand $\overline{u} = (u, \dot{u}, \dots, u^{(k)})$ to
denote a finite but arbitrary number of derivatives of $u$.

\begin{definition}[Definition 3 from \cite{martin2006flat}]
    The system $(\mathfrak{M}, F)$ is \textit{flat} at $p \in \mathfrak{M}$ (resp. \textit{flat}) if it is equivalent at $p$ (resp. equivalent) to a trivial system, i.e., there exists an endogenous transformation $\Psi$ of the form $\Psi(x, u, \dot{u}, \dots) = (h(x, \overline{u}), \dot{h}(x, \overline{u}), \dots)$. Here $h(x, \overline{u})$ is referred to as a \textit{flat output}.
\end{definition}

As a corollary of the following theorem, we have that for a flat system, the input dimension is equal to the dimension of the flat output.

\begin{theorem}[Theorem 1 from \cite{martin2006flat}]
    If two systems $(X \times U \times \mathbb{R}_m^\infty, F)$ and $(Y \times V \times \mathbb{R}_s^\infty, G)$ are equivalent, then they have the same number of inputs, i.e., $m=s$.
\end{theorem}

Equivalence between systems has an alternative characterization in terms of dynamic feedback and coordinate change, as summarized by the following definition and theorem.
\begin{definition}[Definition 5 in \cite{martin2006flat}]
    Consider the dynamics $\dot{x} = f(x, u)$. We say the feedback $u = \kappa(x, z, w), \: \dot{z} = a(x, z, w)$ is \textit{endogenous} if the open-loop dynamics $\dot{x} = f(x, u)$ is equivalent to the closed-loop dynamics $\dot{x} = f(x, \kappa(x, z, w)), \: \dot{z} = a(x, z, w)$.
\end{definition}

\begin{theorem}[Theorem 3 in \cite{martin2006flat}]\label{theorem:equiv-characterization-flat-systems}
    Two dynamics $\dot{x} = f(x, u)$ and $\dot{y} = g(y, v)$ are equivalent if and only if $\dot{x} = f(x, u)$ can be transformed by endogenous feedback and coordinate change into
    \begin{align}
        \dot{y} = g(y, v), \: \dot{v} = v_1, \: \dot{v}_1 = v_2, \: \dots, \: \dot{v}_\mu = w,
    \end{align}
    for some large enough integer $\mu$, and vice versa (possibly extending to a different number of derivatives).
\end{theorem}

In \cite{wardi2024tracking}, the authors assume that $x(t)$ is available at each time $t$. We make the additional assumption that the output of the system which we are interested in regulating to a given reference signal is also the flat output, and that the reference $r$ is a reference for the flat output. In the next section, we make use of the definitions and preliminary results recalled in this section in order to derive conditions for the stability of a tracking control objective applied to a differentially flat system. Moreover, for ease of notation, we will often drop the dependence on time $t$ of input, state, and output signals, unless it is necessary to resolve ambiguities.

\section{Main Results}
\label{sec:main}

In this section, we will provide sufficient conditions for the $\alpha$-stability of the trivial system. Since equivalent trajectories in the flat and real dynamics evolve in parallel as per Fig. \ref{fig:b}, $\alpha$-stability of the flat system follows from the equivalence between the NR controller and the one we propose in this paper. Moreover, we will characterize the cases in which the controller we propose is equivalent to a NR controller.

Given a flat system $\dot{x} = f(x, u)$ with $x \in \mathbb{R}^n, u \in \mathbb{R}^m$, as per Section~\ref{sec:background:flatness}, its trajectories are fully determined by the trajectories of the trivial system with state $\tilde{y} = \begin{bmatrix}y^\top & \dot{y}^\top & \dots & y^{(k)^\top}\end{bmatrix}^\top$ and input $\nu = y^{(k+1)}$ for $y, \nu \in \mathbb{R}^m$ and $k$ finite. Let $h$ denote the endogenous mapping taking $x, \overline{u}$ to the flat output $y$, as per the remark following Theorem \ref{theorem:equiv-characterization-flat-systems}, and $a$ the endogenous mapping taking $y, \dot{y}, \dots, y^{(k+1)}$ to $u$.

We propose the following controller:
\begin{equation}
\begin{aligned}
    \hat{y}(t+T) &= g(\tilde{y}, \nu) \\
    &= g(h(x, \overline{u}), \dot{h}(x, \overline{u}), \dots, h^{(k+1)}(x, \overline{u})) \\
    \dot{\nu} &= y^{(k+2)} = \alpha \left( \frac{\partial g}{\partial \nu} \right)^{-1} e(t+T) \\
    \dot{u} &= \frac{da}{dt}(\tilde{y}, \nu, \dot{\nu}),
\end{aligned}
\label{eq:closed-loop-controller}
\end{equation}
where the predictor is defined as in section \ref{sec:background:NR}, and $\alpha > 1$ is a constant scaling factor. Note that if the dependence on $\nu$ in $a$ is trivial, no new value of $\dot{u}$ is obtained, and so we make the assumption that the input is dynamically extended to nontrivially depend on $\nu$. Refer to Fig. \ref{fig:a} for a diagrammatic description.

Equation~\eqref{eq:closed-loop-controller} encompasses the total of the online computation: once the endogenous transformation is known, both the transformation, its inverse and the computation of the input in the flat dynamics are computable algebraically. It is worthwhile to point out that, in \cite{wardi2024tracking}, robustness to errors in the inverse Jacobian of $g$ in \eqref{eq:closed-loop-controller} is shown for a generalization of the NR controller, which holds for the presented controller.

In section \ref{sec:suff_conds} we derive trivially checkable sufficient conditions for $\alpha$-stability of the trivial dynamics under the application of the controller described by \eqref{eqn:controller}, and in section \ref{subsec:closed-loop-convergence} we show a modification of the above controller for which closed loop convergence can be shown. In section \ref{subsec:NR-on-original-dynamics}, under certain assumptions on the system, we characterize the equivalence between the controller \eqref{eq:closed-loop-controller} we proposed above and the NR controller from \cite{wardi2024tracking}. This equivalence will be illustrated also by means of simulative examples in Sections~\ref{sec:examples} and \ref{sec:simulations}.

\subsection{Sufficient Conditions for $\alpha$-stability of the Trivial System}
\label{sec:suff_conds}

The trivial system is linear and evolves according to the following dynamics:
\begin{equation}
	\begin{cases}
		\dot{\tilde{y}} = A \tilde{y} + B \nu \\
		y = C\tilde{y},
	\end{cases}
    \label{eq:trivial_system_dynamics}
\end{equation}
where
\begin{align}
    A &= \begin{bmatrix}
        0 & I_m & \dots & 0 \\
        \vdots & \vdots & \ddots & \vdots \\
        0 & 0 & \dots & I_m \\
        0 & 0 & \dots & 0 \\
    \end{bmatrix}, \: B = C^\top = \begin{bmatrix}
        0 \\
        \vdots \\
        0 \\
        I_m \\
    \end{bmatrix}.
    \label{eq:trivial_system_dynamicsAB}
\end{align} 
Applying the controller \eqref{eqn:controller} to the dynamics \eqref{eq:trivial_system_dynamics}, we obtain a system of the form \eqref{eqn:augmented_state_evolution}. Substituting the values of $A, B, C$ and computing $P_\alpha(s)$, one obtains (see Appendix):
\begin{align}
    P_{\alpha}(s) = \left( s^{k+2} + \alpha \frac{(k+1)!}{T^{k+1}} \left( \sum_{i=0}^{k+1} \frac{T^i}{i!} s^i \right) \right)^m.
    \label{eq:Palpha}
\end{align}
$P_0(s) = \left( \frac{(k+1)!}{T^{k+1}} \left( \sum_{i=0}^{k+1} \frac{(sT)^i}{i!} \right) \right)^m$ is a finite Maclaurin expansion of $e^{Ts}$ and will have all its roots in the open left half-plane if and only if $k \leq 3$. $Q(s) = \sum_{i=0}^m {m \choose i} s ^ i = (s+1)^m$ by the binomial theorem, so $Q(s)$ always has its roots in the open LHP. Therefore by Theorem \ref{theorem:suff-conds}, the trivial system is $\alpha$-stable when $k \leq 3$. 

\begin{remark} When $k + 2 > 5$ derivatives of the flat output are needed for the inverse endogenous transformation, the following approaches can be employed:
    \begin{enumerate}
        \item $\alpha$-stability can be shown directly from the definition, bypassing the sufficient conditions shown in \cite{wardi2024tracking}
        \item Since the dynamics of the flat outputs are arbitrary, one may consider assigning them an asymptotically stable dynamics, rather than the stable but not asymptotically stable trivial dynamics \eqref{eq:trivial_system_dynamics}--\eqref{eq:trivial_system_dynamicsAB}. The following, with appropriate choices of $K_i$, $i=1,\ldots,k+1$ would suffice:
        \begin{align}
            A &= \begin{bmatrix}
                0 & I_m & \dots & 0 \\
                \vdots & \vdots & \ddots & \vdots \\
                0 & 0 & \dots & I_m \\
                -K_1 I_m & -K_2 I_m & \dots & -K_{k+1} I_m \\
            \end{bmatrix}.
        \end{align}
    \end{enumerate}
\end{remark}

\subsection{Convergence of the closed-loop controller}
\label{subsec:closed-loop-convergence}

The controller described in \eqref{eq:closed-loop-controller} has experimentally been verified to converge in the examples of section \ref{sec:examples}. Although we did not provide a proof of convergence for this controller, we will show a modified controller and the proof of its convergence under the assumption that the flat dynamics are $\alpha$-stable and $r(t)$ is constant.

The modified controller is of the form of Eq. (41) in \cite{wardi2024tracking}, where $\nu$ in the flat dynamics is modified as follows:
\begin{align}
    \dot{\nu} = \left( \frac{\partial g}{\partial \nu} \right) ^{-1} \left(\alpha (r(t+T) - \hat{y}(t+T)) - \left( \frac{\partial g}{\partial \tilde{y}} \dot{\tilde{y}} \right) \right).
\end{align}
We will employ the Lyapunov method in order to show that tracking is achieved. Define the Lyapunov function $V(t) = \frac{1}{2} \|e(t+T)\|^2$, where $e$ is defined as in \eqref{eqn:controller}. The time derivative of $V$ along the trajectories of the system \eqref{eq:trivial_system_dynamics} evaluates to:
\begin{align}
    \dot{V}(t) &= e(t+T)^\top (\dot{r}(t+T) - \dot{\hat{y}}(t+T)) \\
    &= e(t+T)^\top \left( \dot{r}(t+T) - \left(\frac{\partial g}{\partial \tilde{y}} \dot{\tilde{y}} + \frac{\partial g}{\partial \nu} \dot{\nu} \right) \right) \\
    &= e(t+T)^\top \left( \dot{r}(t+T) - \alpha e(t+T) \right).
\end{align}
When $r(t)$ is constant, we have $\dot{V}(t) = -2 \alpha V(t)$, which results in exponential stability of the controlled system \cite{khalil2002nonlinear}.

\subsection{Equivalence with the NR Controller}
\label{subsec:NR-on-original-dynamics}

The controller introduced in the previous section is, under certain conditions, equivalent to the NR controller introduced in \cite{wardi2024tracking}. With the goal of deriving these conditions, let $\Phi$ be the endogenous transformation mapping the flat output dynamics to the ones of the system state, $x, u$. Let $\Psi$ denote the inverse mapping, i.e., the mapping from the system state dynamics to the flat output dynamics. In the following, we will use the following shorthand notation: $z := \begin{bmatrix} x^\top & u^\top \end{bmatrix}^\top$ and $w := \begin{bmatrix} \tilde{y}^\top & \nu^\top \end{bmatrix}^\top$. Our analysis will depend on the following assumptions on the system.

\begin{assumption} \label{a0}
    $\Phi$ does not depend on derivatives of $\nu$, i.e., $\Phi = \Phi(\tilde{y}, \nu)$.
\end{assumption}

\begin{assumption} \label{a1}
    $y=h(x, \bar{u})$ does not depend on derivatives of u, i.e., $h(x, \bar{u}) = h(x, u)$.
\end{assumption}

\begin{assumption} \label{a2}
    $x$ does not depend on $\nu$ in $\Phi$.
\end{assumption}

\begin{assumption} \label{a3}
    $\tilde{y}, x$ have the same dimension, and $\frac{d \Phi}{d w}$ has full rank.
\end{assumption}

The above assumptions are easily verifiable on typical tracking control objectives. Under Assumptions \ref{a0}, \ref{a1}, \ref{a2} and \ref{a3}:
\begin{align}
    \frac{d \Psi}{d z} &= \left( \frac{d \Phi}{d w} \right)^{-1} = {\begin{bmatrix}
        \frac{\partial x}{\partial \tilde{y}} & 0 \\
        \frac{\partial u}{\partial \tilde{y}} & \frac{\partial u}{\partial \nu}
    \end{bmatrix}}^{-1} \\
	&= \begin{bmatrix}
        \left( \frac{\partial x}{\partial \tilde{y}} \right)^{-1} & 0 \\
        - \left( \frac{\partial u}{\partial \nu} \right)^{-1} \frac{\partial u}{\partial \tilde{y}} \left( \frac{\partial x}{\partial \tilde{y}} \right)^{-1} & \left( \frac{\partial u}{\partial \nu} \right)^{-1}
    \end{bmatrix}.
\end{align}
Matching blocks, we see that $\frac{\partial \tilde{y}}{\partial u} = 0$ and $\frac{\partial \nu}{\partial u} = \left( \frac{\partial u}{\partial \nu} \right)^{-1}$. Therefore $\frac{\partial g}{\partial u} = \frac{\partial g}{\partial \nu} \frac{\partial \nu}{\partial u}$ and the time derivative of $u$ can be evaluated as follows:
\begin{align}
    \dot{u} &= \frac{\partial u}{\partial \tilde{y}} \dot{\tilde{y}} + \frac{\partial u}{\partial \nu} \dot{\nu} \\
    &= \frac{\partial u}{\partial \tilde{y}} \dot{\tilde{y}} + \frac{\partial u}{\partial \nu} \alpha \left( \frac{\partial g}{\partial \nu} \right)^{-1} e(t+T) \\
    &= \frac{\partial u}{\partial \tilde{y}} \dot{\tilde{y}} + \alpha \left( \frac{\partial g}{\partial u} \right)^{-1} e(t+T).
\end{align}
The terms in the expression of $\dot{u}$ containing a factor of $\alpha$ are of the form of the NR controller in \eqref{eqn:controller} applied directly to the real dynamics, and any additional terms are of the form $\frac{\partial u}{\partial \tilde{y}} \dot{\tilde{y}}$.

\section{Tracking Control for Vehicle Dynamics} \label{sec:examples}

This section is devoted to the application of the controller developed in the previous section to motion models typically adopted for autonomous mobile platforms, such as ground mobile robots used in warehouse automation \cite{siegwart2011introduction} and autonomous vehicles \cite{rajamani2006vehicle}.

\subsection{Kinematic Unicycle Model}

\begin{figure}
\centering
\includegraphics[width=0.2\textwidth]{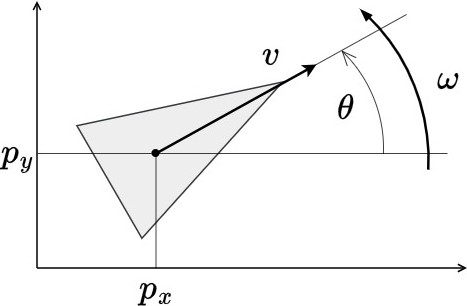}
\caption{Kinematic unicycle: the triangle represents the single wheel of the unicycle, $p_x$ and $p_y$ denote the coordinates of the position of the center of the wheel, and $\theta$ its heading. The inputs $v, \omega$---longitudinal and angular velocities, respectively---are marked in thick black arrows.}
\label{fig:unicycle}
\end{figure}

With reference to Fig.~\ref{fig:unicycle}, the motion of the kinematic unicycle is described by the following system of differential equations:
\begin{align}
	\dot{x} = \begin{bmatrix}
               v \cos{\theta} &
               v \sin{\theta} &
               \omega
            \end{bmatrix}^\top,
	\label{eq:unicycle}
\end{align}
where $x = \begin{bmatrix}
	p_x &
	p_y &
	\theta
\end{bmatrix}^\top$ and $u = \begin{bmatrix}
	v &
	\omega
\end{bmatrix}^\top$. The system is differentially flat \cite{murray1995differential,fuchshumer2005nonlinear}, with flat output $y = \begin{bmatrix} p_x & p_y \end{bmatrix}^\top =: \begin{bmatrix} x_1 & x_2 \end{bmatrix}^\top$. Applying the NR controller to $\nu = \dot{y}$, one obtains the dynamically defined controller for the flat output dynamics:
\begin{align}
   \dot {\nu} &= \left( \frac{\partial g}{\partial \nu} \right)^{-1} (r(t + T) - \hat{y}(t+T)) \\
   &= \frac{\alpha}{T} (r(t + T) - (p(t) + T\dot{p}(t))).
\end{align}

The inputs $v, \omega$ to the dynamics in \eqref{eq:unicycle} can be computed from $y, \nu, \dot{\nu}$ as follows:
\begin{align}
   v &= \sqrt{\dot{p}_x^2 + \dot{p}_y^2} = \sqrt{\nu_1^2 + \nu_2^2} \\
   \omega &= \frac{\ddot{p}_y \dot{p}_x - \dot{p}_y \ddot{p}_x}{\dot{p}_x^2 + \dot{p}_y^2} = \frac{\dot{\nu}_2 \nu_1 - \nu_2 \dot{\nu}_1}{\nu_1^2 + \nu_2^2}.
\end{align}

Since $v$ does not depend on $\dot{\nu}$, we dynamically extend $v$ to $\dot{v}$, which is computed from the parameters of the flat dynamics as $\dot{v} = \frac{\dot{p}_x \ddot{p}_x + \dot{p}_y \ddot{p}_y}{\sqrt{\dot{p}_x^2 + \dot{p}_y^2}} = \frac{\nu_1 \dot{\nu}_1 + \nu_2 \dot{\nu}_2}{\sqrt{\nu_1^2 + \nu_2^2}}$. Then,
\begin{equation}
    \begin{bmatrix}
        \dot{v} \\
        \omega
    \end{bmatrix} = \frac{\alpha}{T} \begin{bmatrix}
        1 & 0 \\
        0 & \frac{1}{v}
    \end{bmatrix} R^{-1}(\theta) e(t+T) + \begin{bmatrix}
        - \alpha v \\
        0
    \end{bmatrix},
\end{equation}
where $R(\theta)$ denotes a clockwise rotation by angle $\theta$.

To satisfy the assumptions in section \ref{subsec:NR-on-original-dynamics}, we take $\dot{v}, \omega$ to be the dynamically extended inputs of the real dynamics, i.e., $u = \begin{bmatrix} \theta & v \end{bmatrix}^\top$. Since only $k + 2 = 2$ derivatives of $y$ are needed for the endogenous transformation, we see that the flat dynamics are $\alpha$-stable.
Assumption \ref{a1} is satisfied for the unicycle, and by taking $u$ as above, Assumptions \ref{a0}, \ref{a2} and \ref{a3} are also satisfied. As a matter of fact, consider the following predictor
\begin{align}
    \hat{y}(t+T) &= \begin{bmatrix}
        p_x + Tv \cos \theta \\
        p_y + Tv \sin \theta,
    \end{bmatrix}
\end{align}
the NR controller applied directly to the dynamics of the unicycle is of the same form as the one obtained by applying the controller described in \eqref{eq:closed-loop-controller}. Note that $\frac{\partial u}{\partial \tilde{y}} = 0$, so this is consistent with the results described in section \ref{subsec:NR-on-original-dynamics}.

\subsection{Dynamic Bicycle Model}

\begin{figure}
\centering
\includegraphics[width=0.32\textwidth]{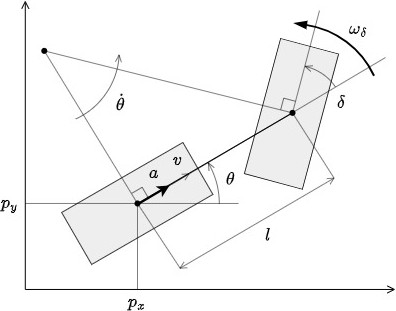}
\caption{Dynamic bicycle: the two wheels are represented by gray-shaded rectangles, joined by the wheelbase. $p_x$, $p_y$, $\theta$, and $v$ denote the coordinates of the rear axle, the heading, and the velocity of the back wheel, respectively. $\delta$ represents the steering angle, and $l$ is the length of the wheelbase. The inputs $a, \omega_\delta$---the acceleration of the rear wheel and the steering speed of the front wheel, respectively---are marked in thick black arrows.}
\label{fig:bicycle}
\end{figure}

A vehicle dynamics model is considered to be \textit{dynamic} if it accounts for drift, regardless of whether it takes into account tire forces. Here we consider the dynamic bicycle model, depicted in Fig.~\ref{fig:bicycle} and described by the following differential equations: 
\begin{align}
	\dot{x} = \begin{bmatrix}
            v \cos{\theta} &
            v \sin{\theta} &
            \frac{v}{l} \tan{\delta} &
            a &
            \omega_{\delta}
         \end{bmatrix}^\top,
	\label{eq:bicycle}
\end{align}
where $x = \begin{bmatrix}
	p_x &
	p_y &
	\theta &
	v &
	\delta
\end{bmatrix}^\top$ and $u = \begin{bmatrix}
	a &
	\omega_{\delta}
\end{bmatrix}^\top$. These dynamics are differentially flat with flat output $y = \begin{bmatrix} p_x & p_y \end{bmatrix}^\top =: \begin{bmatrix} x_1 & x_2 \end{bmatrix}^\top$. Derivatives of $y$ up to the third order are necessary to compute the input $\omega_\delta$, so we let the state $\tilde{y} = \begin{bmatrix} y^\top & \dot{y}^\top \end{bmatrix}^\top$ and the input $\nu = \ddot{y}$. The NR controller is applied to compute $\dot{\nu}$ and the following expression for the dynamically defined controller of the flat output dynamics is obtained:
\begin{align}
   \dot{\nu} &= \left( \frac{\partial g}{\partial \nu} \right)^{-1} (r(t + T) - \hat{y}(t+T)) \\
   &= \frac{2 \alpha}{T^2} \left( r(t + T) - \left(p(t) + T\dot{p}(t) + \frac{T^2}{2} \ddot{p}(t) \right) \right).
\end{align}

The longitudinal acceleration $a$ can be calculated from $\tilde{y}$, $\nu$ as follows:
\begin{align}
    a = \frac{d}{dt} \sqrt{\dot{p}_x^2 + \dot{p}_y^2}
    = \frac{\dot{p}_x \ddot{p}_x + \dot{p}_y \ddot{p}_y}{\sqrt{\dot{p}_x^2 + \dot{p}_y^2}}
    = \frac{\dot{y_1} \nu_1 + \dot{y_2} \nu_2}{\sqrt{\dot{y_1}^2 + \dot{y_2}^2}}.
\end{align}
For ease of notation, let $q := \ddot{p}_y \dot{p}_x - \ddot{p}_x \dot{p}_y$. Its time derivative is $\dot{q} = \dddot{p}_y \dot{p}_x - \dddot{p}_x \dot{p}_y$. Recall that $v = \sqrt{\dot{p}_x^2 + \dot{p}_y^2}$ and $av = \dot{p}_x \ddot{p}_x + \dot{p}_y \ddot{p}_y$. Then,
\begin{align}
    \omega_{\delta} = \frac{d}{dt} \tan^{-1} \left( \frac{l}{v} \dot{\theta} \right)
    = \frac{l v}{v^6 + l^2 q^2} ( \dot{q} v^2 - 3 q a v ).
\end{align}

Since $a$ does not depend on $\dot{\nu}$, we dynamically extend the input to $\dot{a}$. Note that $\ddot{p}_x = \frac{d}{dt} v \cos \theta = a \cos \theta - \frac{v^2}{l} \sin \theta \tan \delta$ and $\ddot{p}_y = \frac{d}{dt} v \sin \theta = a \sin \theta + \frac{v^2}{l} \cos \theta \tan \delta$. Let $R(\theta)$ denote a clockwise rotation by angle $\theta$. Then,
\begin{align}
    \begin{bmatrix}
        \dot{a} \\
        \omega_\delta
    \end{bmatrix} &= \frac{2 \alpha}{T^2} \Bigg( \!\! \Bigg( \! \begin{bmatrix}
        1 & 0 \\
        0 & \frac{l \cos^2\delta}{v^2} \\
    \end{bmatrix} R^{-1}(\theta) e(t+T) \!\Bigg) \!-\! \begin{bmatrix}
        Tv + \frac{T^2}{2}a \\
        \sin \delta \cos \delta
    \end{bmatrix} \!\Bigg) \\
    &+ \begin{bmatrix}
        \frac{v^3}{l} \tan^2 \delta \\
        - 3 \frac{a}{v} \cos \delta \sin \delta
    \end{bmatrix}.
    \label{eq:dynbicycleNRflat}
\end{align}

To satisfy the assumptions in section \ref{subsec:NR-on-original-dynamics}, we take $\dot{a}, \omega_\delta$ to be the dynamically extended inputs of the real dynamics, i.e., $u = \begin{bmatrix} \delta & a \end{bmatrix}^\top$. Since only $k + 2 = 3$ derivatives of $y$ are needed for the endogenous transformation, we see that the flat dynamics are $\alpha$-stable.

We see that Assumption \ref{a1} is satisfied for the dynamic bicycle, and by taking $u$ as above, Assumptions \ref{a0}, \ref{a2} and \ref{a3} are satisfied. Using the following predictor
\begin{align}
    \hat{y}(t+T) &= \begin{bmatrix}
        p_x + Tv \cos \theta + \frac{T^2}{2} \left( a \cos \theta - \frac{v^2}{l} \sin \theta \tan \delta \right) \\
        p_y + Tv \sin \theta + \frac{T^2}{2} \left( a \sin \theta + \frac{v^2}{l} \cos \theta \tan \delta \right)
    \end{bmatrix},
\end{align}
the NR controller applied to the bicycle dynamics \eqref{eq:bicycle} matches the form of the controller defined in \eqref{eq:closed-loop-controller} (see terms multiplied by $\alpha$ in \eqref{eq:dynbicycleNRflat}). Moreover, it is easy to check that the additional terms are of the form $\frac{\partial u}{\partial \tilde{y}} \dot{\tilde{y}}$.

\section{Simulation Results}
\label{sec:simulations}

The kinematic unicycle and dynamic bicycle models were simulated using the controller proposed in this paper to follow desired reference trajectories, namely a sine wave of the form $r(t) = \begin{bmatrix} \frac{1}{5} t & 10 \sin \left(\frac{2 \pi}{50} t \right) \end{bmatrix}^\top$ and a spiral of the form $r(t) = e^{0.0125s} \begin{bmatrix} \cos(0.25s) & \sin(0.25s) \end{bmatrix}^\top, s = 284 - t$, both over the time range $0 \leq t \leq 100$, for $t$ in seconds. We apply a time varying reference here, in contrast to the constant reference assumed in section \ref{subsec:closed-loop-convergence} to illustrate how stability guarantees can be extended to the case of a non-constant reference.

Distances are measured in meters, and the simulation is updated using forward Euler with a time step of $0.001$s. In all simulations, the vehicle (unicycle or bicycle) starts at a heading $\theta = \frac{\pi}{2}$. In the simulation following the sine wave reference, the vehicle starts at position $(12, -4)$, while for the simulation with the spiral reference, the vehicle starts at $(-12, -13)$. The other parameters required to be able to reproduce the simulations have been set to the following values. The length of the bicycle wheelbase is $l = 2$, the controller parameters are $\alpha = 30$ and $T = 0.8$. For the unicycle, we set $\alpha = 100$ and $T = 0.02$.

\begin{figure}
	\centering
	\subfloat[]{{\includegraphics[height=0.42\linewidth]{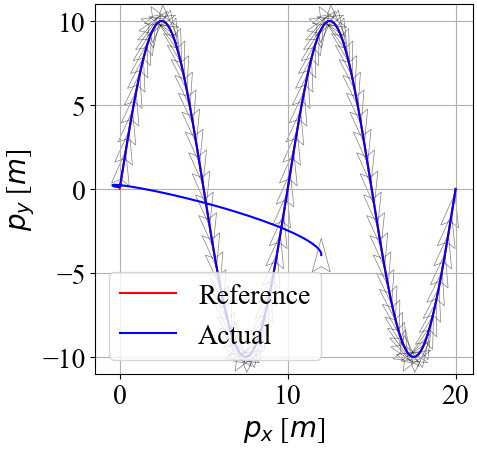} }}%
	\quad
	\subfloat[]{{\includegraphics[height=0.42\linewidth]{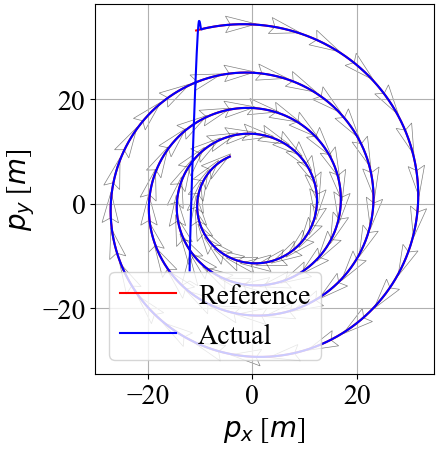} }}%
	\caption{Kinematic unicycle simulated to follow two reference trajectories (sine wave on the left and spiral on the right). The reference trajectories are plotted in red, while the unicycle trajectory is in blue (overlapping with red). Arrows show the heading of the unicycle at 1 second time intervals.}%
	\label{fig:unicycle_trajectory}%
\end{figure}
\begin{figure}
	\centering
	\includegraphics[width=0.45\textwidth]{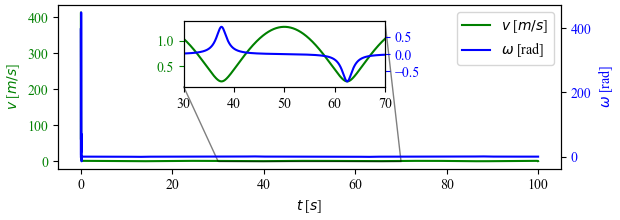}
	\caption{Control inputs to the unicycle computed using the proposed controller to track a sine wave reference trajectory.}
	\label{fig:unicycle_sine_wave_inputs}
\end{figure}
\begin{figure}
	\centering
	\includegraphics[width=0.45\textwidth]{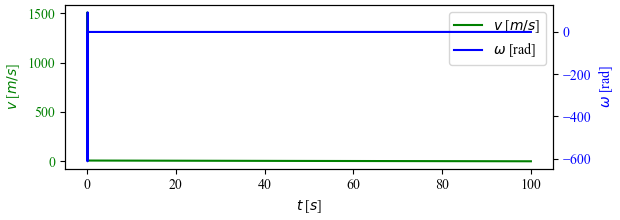}
	\caption{Control inputs to the unicycle computed using the proposed controller to track a spiral reference trajectory.}
	\label{fig:unicycle_spiral_inputs}
\end{figure}

\begin{figure}
	\centering
	\subfloat[\centering label 1]{{\includegraphics[height=0.42\linewidth]{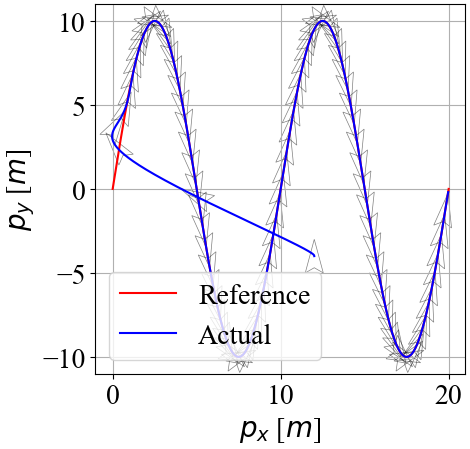} }}%
	\quad
	\subfloat[\centering label 2]{{\includegraphics[height=0.42\linewidth]{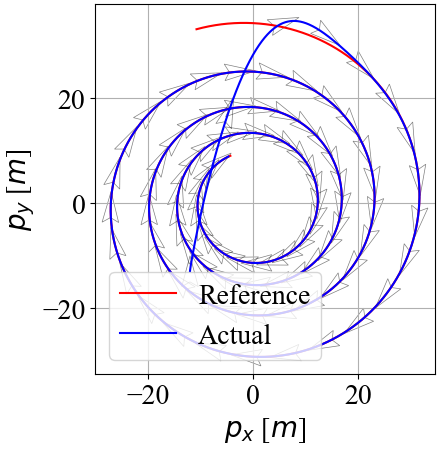} }}%
	\caption{Dynamic bicycle simulated to follow two reference trajectories (sine wave on the left and spiral on the right). The reference trajectories are plotted in red, while the trajectory of the rear axle of the bicycle are in blue. Arrows show the heading of the bicycle at 1 second time intervals.}%
	\label{fig:bicycle_trajectory}%
\end{figure}
\begin{figure}
	\centering
	\includegraphics[width=0.45\textwidth]{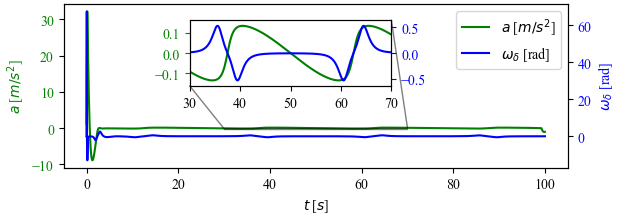}
	\caption{Control inputs to the bicycle computed using the proposed controller to track a sine wave reference trajectory.}
	\label{fig:bicycle_sine_wave_inputs}
\end{figure}
\begin{figure}
	\centering
	\includegraphics[width=0.45\textwidth]{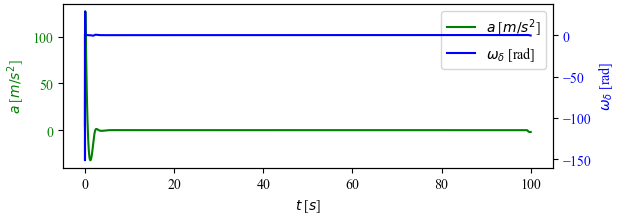}
	\caption{Control inputs to the bicycle computed using the proposed controller to track a spiral reference trajectory.}
	\label{fig:bicycle_spiral_inputs}
\end{figure}
\begin{figure}
	\centering
	\includegraphics[width=0.45\textwidth]{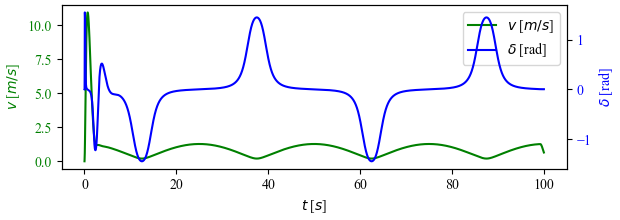}
	\caption{Velocity and heading of bicycle while tracking the sine wave reference trajectory.}
	\label{fig:bicycle_sine_wave_extra}
\end{figure}
\begin{figure}
	\centering
	\includegraphics[width=0.45\textwidth]{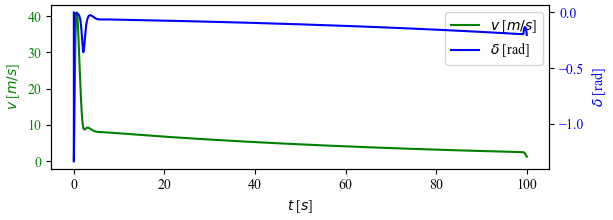}
	\caption{Velocity and heading of bicycle while tracking the spiral reference trajectory.}
	\label{fig:bicycle_spiral_extra}
\end{figure}

The results of the unicycle simulations can be seen in Figures \ref{fig:unicycle_trajectory}--\ref{fig:unicycle_spiral_inputs}. The results of the bicycle simulation can be seen in Figure \ref{fig:bicycle_trajectory}--\ref{fig:bicycle_spiral_extra}. Both vehicle models converge rapidly to the reference in simulation, with small overshoot. The large inputs observed at the beginning are due to the large initial magnitude of error. This effect can be eliminated by planing a path that originates from the initial state of the vehicle. Note that in the simulation of the bicycle, $\delta$ was not clamped, resulting in values of the steering angle that may not be physically realizable. Input constraints have not been considered in this work, however approaches similar to the ones in \cite{ames2020integral} can be employed to account for them.

\section{Conclusions}
\label{sec:conclusions}

In this paper, we presented a computationally efficient, rapidly converging, controller suitable for the tracking control of differentially flat systems. We derived sufficient conditions for the $\alpha$-stability of the dynamics of the flat outputs of these systems and we presented a proof of convergence for a modified controller not further analyzed in this paper. Further, we presented a set of conditions under which this controller we propose is equivalent to applying the Newton-Raphson tracking controller directly to the nonlinear dynamics of the state. Simulations employing kinematic and dynamic models of mobile robots illustrate the properties of the defined controller.

\appendix

\begin{align}
	(sI-A)^{-1} &= \begin{bmatrix}
		\frac{1}{s} I & \frac{1}{s^2} I & \dots & \frac{1}{s^{k+1}} I \\
		0 & \frac{1}{s} I & \dots & \frac{1}{s^k} I \\
		\vdots & \ddots & \ddots & \vdots \\
		0 & \dots & 0 & \frac{1}{s} I \\
	\end{bmatrix},\\
	e^{A \tau} &= \begin{bmatrix}
		I & \tau I & \dots & \frac{\tau^k}{k!} I \\
		0 & I & \dots & \frac{\tau^{k-1}}{(k-1)!} I \\
		\vdots & \ddots & \ddots & \vdots \\
		0 & \dots & 0 & I \\
	\end{bmatrix},\\
	C e^{AT} (sI-A)^{-1} B &= \left( \sum_{i=0}^k \frac{T^i}{i!} \frac{1}{s^{k + 1 - i}} \right) I,\\
	\left( C \int_0^T e^{A \tau} d \tau B \right)^{-1}
	&= \frac{(k+1)!}{T^{k+1}} I.
\end{align}

With the expressions above, the expression of $P_{\alpha}(s)$ in \eqref{eq:Palpha} is obtained as follows:
\begin{align}
	P_{\alpha}(s) &= \det  \begin{bmatrix}
		sI-A & -B \\
		\alpha \left( C \int_0^T e^{A \tau} d \tau B \right)^{-1} C e^{AT} & (s + \alpha) I \\
	\end{bmatrix}\\
	&= \det(sI-A) \det \Bigg( (s + \alpha) I \\
	&\quad+ \alpha \left( C \int_0^T e^{A \tau} d \tau B \right)^{-1} C e^{AT} (sI-A)^{-1} B \Bigg)\\
	&= s^{m(k+1)} \det \bigg( (s+\alpha) I \\
	&\quad+ \alpha \frac{(k+1)!}{T^{k+1}} \left( \sum_{i=0}^k \frac{T^i}{i!} \frac{1}{s^{k + 1 - i}} \right) I \bigg)\\
	&= \left( s^{k+2} + \alpha \frac{(k+1)!}{T^{k+1}} \left( \sum_{i=0}^{k+1} \frac{T^i}{i!} s^i \right) \right)^m.
\end{align}

\bibliographystyle{IEEEtran}
\bibliography{bib/references}

\end{document}